\documentclass{article}
\usepackage{amssymb}
\usepackage{amsfonts}
\usepackage{amsmath}

\newtheorem{theorem}{Theorem}

\newtheorem{corollary}[theorem]{Corollary}

\newtheorem{definition}[theorem]{Definition}

\newtheorem{lemma}[theorem]{Lemma}

\newtheorem{proposition}[theorem]{Proposition}

\input{tcilatex}
\begin{document}

\title{Tensor calculus and deformation theory on \\
on a local Lie group}
\author{Erc\"{u}ment Orta\c{c}gil}
\maketitle

\begin{abstract}
Historically tensor calculus emerged in an attempt to formalize Riemann's
ideas. We show that tensor calculus can be based also on Lie's idea of a
transformation group and this approach leads quite naturally to the concept
of deformation of a transformation group and the Kodaira-Spencer map.
\end{abstract}

\section{Introduction}

The study of the deformations of complex structures for higher dimensional
complex manifolds is initiated by K.Kodaira and D.C.Spencer in their
foundational work [7]. The book [6] is an excellent source for learning
deformation theory from one of its creators.

A complex structure is a particular pseudogroup. The ideas and methods in
the above works are later generalized by D.C.Spencer to the deformations of
transitive pseudogroups in [12]. However, it seems that these later works
(which are quite difficult to read) did not have much impact on the
mathematical community. Deformations of transitive pseudogroups are taken up
again in [11] from a different viewpoint. This work also did not attract
much attention possibly due to the attitude of this author in his later
works. Ironically, deformation theory flourished in the hands of algebraists
and algebraic geometers with the impetus from A.Grothendieck and this
situation prevails today.

Finite order transitive pseudogroups are the simplest types of pseudogroups
and they are equivalent to transitive Lie group actions, that is,
homogeneous spaces. Therefore, a complex manifold belongs to this category
if and only if it is a homogeneous complex manifold ([13]). In spite of
their simplicity, finite order pseudogroups are so rich that they can be
taken as the basis of the generalization of Klein's Erlangen Program once
they are decoupled with the concept of curvature ([10]). The simplest finite
order pseudogroups are, of course, those whose order is smallest, that is,
local Lie groups. It is shown in [9] that a local Lie group is not
necessarily contained in a global Lie group. In fact, just the opposite is
true: a Lie group is a special (globalizable, see Section 2) local Lie group
([1]). How does the deformation theory of a local Lie group look like? This
question has been the main motivation for the present author who, as a
newcomer to this field, is deeply impressed by the personal and pedagogical
aspects of the remarkable book [6]. However, it was quite surprising for him
to realize that to study this problem one must first develop the formalism
of right-left covariant differentitions and curvatures as far as the fine
resolutions of right-left Lie algebras. We believe that this "tensor
calculus" on a local Lie group which we barely touch in Sections 3, 4
neccesitates our reexamination the concept of torsion in modern differential
geometry.

In the following section we outline the theory of local Lie groups initiated
in [9], [1]. Sections 3, 4 give an introduction to the tensor calculus on a
local Lie group. Section 5 introduces the concept of deformation via gauge
transformations and in Section 6 we construct the Kodaira-Spencer map.

\section{Local Lie groups}

In this section we briefly recall local Lie groups. We will continue to use
the notation of [1] and refer to [1] for more details for certain points.

Let $M$ be a smooth manifold with $\dim M\geq 2$ and $\mathcal{U}_{k}$ be
the universal groupoid of order $k$ on $M.$ The elements of $\mathcal{U}_{k}$
are the $k$-jets of local diffeomorphisms of $M.$ We call an element of $%
\mathcal{U}_{k}$ with source at $p$ and target at $q$ a $k$-arrow (from $p$
to $q).$ Therefore $\mathcal{U}_{0}$ is the pair groupoid $M\times M.$ The
relevant universal groupoids for the study of local Lie groups are $\mathcal{%
U}_{0}$ and $\mathcal{U}_{1}.$ The projection homomorphism $\pi :\mathcal{U}%
_{1}\rightarrow \mathcal{U}_{0}$ of groupoids maps the $1$-arrows from $p$
to $q$ to the pair $(p,q).$ A splitting $\varepsilon :$ $\mathcal{U}%
_{0}\rightarrow \mathcal{U}_{1}$ is a homomorphism of groupoids satisfying $%
\pi \circ \varepsilon =id_{\mathcal{U}_{0}}.$ Thus $\varepsilon $ assigns to
any pair $(p,q)$ a $1$-arrow from $p$ to $q$ and this assignment preserves
the composition and inversion of arrows. It is easily checked that $\pi :%
\mathcal{U}_{1}\rightarrow \mathcal{U}_{0}$ admits a splitting if and only
if $M$ is parallelizable. If $p\in (U,x^{i})$ has coordinates $\overline{x}%
^{i}$ and $q\in (V,y^{i})$ has coordinates $\overline{y}^{i},$ then $%
\varepsilon (p,q)$ has the coordinate representation $\varepsilon _{j}^{i}(%
\overline{x}^{1},...,\overline{x}^{n},\overline{y}^{1},...,\overline{y}%
^{n})=\varepsilon _{j}^{i}(\overline{x},\overline{y}),$ $1\leq i,j\leq
n=\dim M.$ Thus we have $\varepsilon _{a}^{i}(y,z)\varepsilon
_{j}^{a}(x,y)=\varepsilon _{j}^{i}(x,z),$ $\varepsilon _{j}^{i}(x,x)=\delta
_{j}^{i}$ and $\varepsilon _{a}^{i}(y,x)\varepsilon _{j}^{a}(x,y)=\delta
_{j}^{i}.$ The curvature $\mathcal{R}(\varepsilon )$ of the splitting $%
\varepsilon $ is defined as the integrability conditions of the $PDE$

\begin{equation}
\frac{\partial f^{i}(x)}{\partial x^{j}}=\varepsilon _{j}^{i}(x,f(x))
\end{equation}

Therefore $\mathcal{R}(\varepsilon )=0$ if and only if the subgroupoid $%
\varepsilon (\mathcal{U}_{0})\subset \mathcal{U}_{1}$ integrates to a
pseudogroup $\mathcal{G}$. In this case, if the local diffeomorphism $%
f:(U,x^{i})\rightarrow (V,y^{i})$ with $f(\overline{x})=\overline{y}$
belongs to $\mathcal{G}$ and $(U,x^{i})$ is simply connected, then $f$ is
uniquely determined on its domain $(U,x^{i})$ by its $1$-arrow $\frac{%
\partial f^{i}(\overline{x})}{\partial x^{j}}$ because $f$ is the unique
solution of the first order nonlinear system (1) with the initial condition $%
\varepsilon (\overline{x},\overline{y}).$ A parallelizable manifold $%
(M,\varepsilon )$ is called a local Lie group if $\mathcal{R}(\varepsilon )=0
$. Sometimes (for instance if $(M,\varepsilon )$ is complete and $M$ is
simply connected) all $f\in \mathcal{G}$ extend (necessarily uniquely) to
global diffeomorphisms of $M$. In this case $\mathcal{G}$ becomes a
transformation group of $M$ which acts simply transitively$,$ that is, a Lie
group. If this happens, $\mathcal{G}$ is called globalizable ([9], [1]).
Therefore, a Lie group is a globalizable local Lie group. If $(M,\varepsilon
)$ is complete, the pseudogroup $\mathcal{G}$ lifts to the universal
covering space of $M$ and globalizes to a Lie group there.

It is very useful to linearize (1). This amounts to passing from the
splitting $\varepsilon :\mathcal{U}_{0}\rightarrow \mathcal{U}_{1}$ to its
infinitesimal splittings $\widetilde{\Gamma },$ $\widehat{\Gamma }%
:T\rightarrow J_{1}T,$ that is, passing from the groupoid $\varepsilon (%
\mathcal{U}_{0})\subset \mathcal{U}_{1}$ to its two algebroids $\widetilde{%
\Gamma }(T),$ $\widehat{\Gamma }(T)\subset J_{1}T$ (note that an abstract
Lie groupoid has only one algebroid and it does not have an adjoint
representation (see [8], [3]). The main point here is that $\varepsilon (%
\mathcal{U}_{0})$ is a very special groupoid). To do this, we define the
local components $\Gamma _{kj}^{i}(x)$ by

\begin{equation}
\Gamma _{kj}^{i}(x)\overset{def}{=}\left[ \frac{\partial \varepsilon
_{j}^{i}(x,y)}{\partial y^{k}}\right] _{y=x}
\end{equation}

The components (2) define two splittings $\widetilde{\Gamma },\widehat{%
\Gamma }$ of the projection $\pi :J_{1}T\rightarrow T$ by

\begin{equation}
\begin{array}{ccc}
1)\text{ }\widetilde{\Gamma }:T\rightarrow J_{1}T\text{ } &  & 2)\text{ \ }%
\widehat{\Gamma }:T\rightarrow J_{1}T \\
(\xi ^{i})\rightarrow (\xi ^{i},\Gamma _{ja}^{i}\xi ^{a}) &  & (\xi
^{i})\rightarrow (\xi ^{i},\Gamma _{aj}^{i}\xi ^{a})%
\end{array}%
\end{equation}

$1)$ and $2)$ give two linear systems of $PDE$'s on the space of vector
fields $\mathfrak{X}(M)$ defined by

\begin{equation}
1^{\prime })\text{ }\widetilde{\nabla }_{j}\xi ^{i}\overset{def}{=}\frac{%
\partial \xi ^{i}}{\partial x^{j}}-\Gamma _{ja}^{i}\xi ^{a}=0\text{ \ \ \ \
\ \ \ \ }2^{\prime })\text{ }\widehat{\nabla }_{j}\xi ^{i}\overset{def}{=}%
\frac{\partial \xi ^{i}}{\partial x^{j}}-\Gamma _{aj}^{i}\xi ^{a}=0\text{ }
\end{equation}

with integrability conditions

\begin{equation}
1^{\prime \prime })\text{ }\widetilde{\mathfrak{R}}_{rj,k}^{i}\overset{def}{=%
}\left[ \frac{\partial \Gamma _{jk}^{i}}{\partial x^{r}}+\Gamma
_{rk}^{a}\Gamma _{ja}^{i}\right] _{[rj]}=0\text{ \ \ \ }2^{\prime \prime })%
\text{ }\widehat{\mathfrak{R}}_{rj,k}^{i}\overset{def}{=}\left[ \frac{%
\partial \Gamma _{kj}^{i}}{\partial x^{r}}+\Gamma _{kr}^{a}\Gamma _{aj}^{i}%
\right] _{[rj]}=0
\end{equation}

\ If $2^{\prime \prime })$ is satisfied, then $2^{\prime })$ is locally
solvable uniquely with the initial condition $(\xi ^{i}(p))$ which can be
assigned arbitrarily at any point $p\in M.$ The same statement applies to $%
1^{\prime \prime })$ and $1^{\prime }).$ In fact, let $\mathfrak{X}%
_{\varepsilon }(M)$ denote the $\varepsilon $-invariant vector fields on $M$
which is a vector subspace of the infinite dimensional Lie algebra $%
\mathfrak{X}(M).$ For a parallelizable manifold $(M,\varepsilon ),$ $%
1^{\prime \prime })$ is always satisfied and a solution of $1^{\prime })$ is
the restriction of some element of $\mathfrak{X}_{\varepsilon }(M)$ to $%
(U,x^{i}).$ Therefore, the globally defined $PDE$ $1^{\prime })$ admits $%
\mathfrak{X}_{\varepsilon }(M)$ as global solutions. Henceforth we will use
the more common sheaf notation $\widetilde{\Theta }$ for $\mathfrak{X}%
_{\varepsilon }(M).$ So $\widetilde{\Theta }$ is a $n$-dimensional vector
space whose elements are \textit{global }vector fields on $M$. The problem
is that $\widetilde{\Theta }$ is not closed under bracket and therefore is
not a Lie algebra. Now a straightforward computation shows that the bracket
of two solutions of $1^{\prime })$ is again a solution of $1^{\prime })$ if
and only if $2^{\prime \prime })$ holds (this computation uses also $%
1^{\prime \prime }))$. Therefore, $\widetilde{\Theta }$ is a Lie algebra if
and only if $\widehat{\mathfrak{R}}=0$ (see (12) and the paragraph following
it)$.$ If $\widehat{\mathfrak{R}}=0,$ the globally defined $PDE$ $2^{\prime
})$ is only locally solvable. These local solutions are also closed under
bracket and form a Lie algebra isomorphic to $\widetilde{\Theta }$ but may
not piece together to give global solutions of $2^{\prime })$ if $M$ is not
simply connected. We denote the sheaf of local solutions of $2^{\prime })$
by $\widehat{\Theta }.$

Now suppose $\widehat{\mathfrak{R}}=0$, $p\in M$ and $\xi ,$ $\eta \in
\widetilde{\Theta }.$ Using (3) we compute

\begin{eqnarray}
\lbrack \xi ,\eta ]^{i}(p) &=&\left[ \xi ^{a}(x)\frac{\partial \eta ^{i}(x)}{%
\partial x^{a}}-\eta ^{a}(x)\frac{\partial \xi ^{i}(x)}{\partial x^{a}}%
\right] _{x=p}  \notag \\
&=&\left[ \xi ^{a}(x)\Gamma _{ab}^{i}(x)\eta ^{b}(x)-\eta ^{a}(x)\Gamma
_{ab}^{i}(x)\xi ^{b}(x)\right] _{x=p}  \notag \\
&=&\xi ^{a}(p)\Gamma _{ab}^{i}(p)\eta ^{b}(p)-\eta ^{b}(p)\Gamma
_{ba}^{i}(p)\xi ^{a}(p)  \notag \\
&=&(\Gamma _{ab}^{i}(p)-\Gamma _{ba}^{i}(p))\xi ^{a}(p)\eta ^{b}(p)  \notag
\\
&=&T_{ab}^{i}(p)\xi ^{a}(p)\eta ^{b}(p)
\end{eqnarray}%
where

\begin{equation}
T_{jk}^{i}\overset{def}{=}\Gamma _{jk}^{i}-\Gamma _{kj}^{i}
\end{equation}

We define a bilinear and alternating operation $[\cdot ,\cdot ]_{p}$ on the
tangent space $T_{p}$ as follows: if $\xi _{p},\eta _{p}\in T_{p}$, then $%
[\xi _{p},\eta _{p}]_{p}$ is defined by the last formula in (6), that is, $%
[\xi _{p},\eta _{p}]_{p}\overset{def}{=}T_{ab}^{i}(p)(\xi _{p})^{a}(\eta
_{p})^{b}.$ The computation in (6) shows that $[\cdot ,\cdot ]_{p}$ turns $%
T_{p}$ into a Lie algebra and the evaluation map $i_{p}:$ $\widetilde{\Theta
}\rightarrow T_{p}$ defined by $\xi \rightarrow \xi (p)$ is a homomorphism
of Lie algebras. Now $i_{p}$ is surjective since we can define the initial
condition arbitrarily and is also injective since the solution is unique.
Thus $i_{p}$ is an isomorphism. We will denote the Lie algebra $%
(T_{p},[\cdot ,\cdot ]_{p})$ by $\widetilde{\Theta }_{p}.$ In short, the Lie
algebra of global vector fields $\widetilde{\Theta }$ can be localized at
any point $p\in M$ and is isomorphic to its localization $\widetilde{\Theta }%
_{p}.$

It s a fundamental fact that $\mathcal{R}(\varepsilon )=0\Leftrightarrow
\widehat{\mathfrak{R}}=0.$ The nontrivial implication $\widehat{\mathfrak{R}}%
=0\Rightarrow \mathcal{R}(\varepsilon )=0$ is the \textit{global }version of
the \textit{traditional }form of the Lie's third fundamental theorem $(LTFT)$
and is quite different in essence from the Cartan's version of $LTFT$ which
is commonly accepted today. For the local Lie group $(M,\varepsilon )$, $%
\widehat{\Theta }$ integrates to a locally transitive pseudogroup whose
local diffeomorphisms belong the globally transitive pseudogroup $\mathcal{G}
$ explained above. Similarly, the Lie algebra $\widetilde{\Theta }$
integrates to another locally transitive pseudogroup $\mathcal{H}$ whose $1$%
-arrows are computed in [1] (see (42) in [1]). If we agree to call the
elements of $\widetilde{\Theta }$ left invariant, then $\mathcal{H}$
corresponds to the right action. According to this convention, elements of $%
\widehat{\Theta }$ become right invariant and their integrals $\mathcal{G}$
corresponds to the left action. However, there is no canonical choice of
left and right for a local Lie group even if it is globalizable in contrast
to abstract Lie groups.

\section{$\protect\varepsilon $-invariance}

Recall that $\widetilde{\Theta }$ denotes $\varepsilon $-invariant vector
fields on $M.$ More generally, we define $\varepsilon $-invariant tensor
fields on a parallelizable manifold $(M,\varepsilon )$ as follows. Since $%
\varepsilon (p,q)$ defines an isomorphism $\varepsilon (p,q)_{\ast
}:T_{p}\rightarrow T_{q}$, $\varepsilon (p,q)_{\ast }$ extends to an
isomorphism on tensor spaces (using the same notation) $\varepsilon
(p,q)_{\ast }:(T_{p})_{r}^{s}\rightarrow (T_{q})_{r}^{s}.$ A tensor field $%
\xi \in T_{r}^{s}$ (as in [1], we use the same notation $E$ for both the
total space and for the space of sections of a bundle $E\rightarrow M)$ is $%
\varepsilon $-invariant if $\varepsilon (p,q)_{\ast }\xi (p)=\xi (q)$, $%
p,q\in M.$ Equivalently, we may fix some $p\in M$ and require $\varepsilon
(p,x)_{\ast }\xi (p)=\xi (x),$ $x\in M.$ This condition is independent of $p$
and is equivalent to the first since $\varepsilon (p,q)=\varepsilon
(x,q)\circ \varepsilon (p,x)$ where $\circ $ denotes composition of $1$%
-arrows. Clearly, a tensor $\xi _{p}\in (T_{p})_{r}^{s}$ extends uniquely to
an $\varepsilon $-invariant tensor field $\xi \in T_{r}^{s}$ with $\xi
(p)=\xi _{p}.$ We call $\xi $ the extension of $\xi _{p}.$

In coordines, if $\xi =(\xi _{rs...t}^{ij...k}(x)),$ then the $\varepsilon $%
-invariance of $\xi $ is expressed by the formula

\begin{equation}
\xi _{rs...t}^{ij...k}(x)=\varepsilon _{a}^{i}(p,x)\varepsilon
_{b}^{j}(p,x)...\varepsilon _{c}^{k}(p,x)\xi
_{de...f}^{ab...c}(x)\varepsilon _{r}^{d}(x,p)\varepsilon
_{s}^{e}(x,p)...\varepsilon _{t}^{f}(x,p)
\end{equation}

Now we define the $\widetilde{\cdot }$ -covariant derivative of a tensor
field by the formula

\begin{eqnarray}
&&\widetilde{\nabla }_{l}\xi _{rs...t}^{ij...k}\text{ }\overset{def}{=}\text{
}\frac{\partial \xi _{rs...t}^{ij...k}}{\partial x^{l}}-\Gamma _{la}^{i}\xi
_{rs...t}^{aj...k}-\Gamma _{la}^{j}\xi _{rs...t}^{ia...k}...-\Gamma
_{la}^{k}\xi _{rs...t}^{ij...a}  \notag \\
&&+\Gamma _{lr}^{a}\xi _{as...t}^{ij...k}+\Gamma _{ls}^{a}\xi
_{ra...t}^{ij...k}...+\Gamma _{lt}^{a}\xi _{rs...a}^{ij...k}
\end{eqnarray}

Our sign convention in (9) is the opposite of the one in tensor calculus
since we define $\Gamma _{kj}^{i}$ by (2) rather than $\left[ \frac{\partial
\varepsilon _{j}^{i}(y,x)}{\partial y^{k}}\right] _{y=x}=-$ $\Gamma
_{kj}^{i} $ but this point is not much important. However, unlike in
Riemannian geometry we should be careful with order of the lower indices of $%
\Gamma _{kj}^{i}$ in (9): the differentiation index $l$ appears always as
the first lower index in $\Gamma _{l\cdot }^{\cdot }$ as in $1^{\prime }).$
Clearly (9) defines a differential operator $\widetilde{\nabla }%
_{X}:T_{r}^{s}\rightarrow T_{r}^{s}$ for all $X\in \mathfrak{X}(M)$ with the
well known properties. Thus we also have the operator $\widetilde{\nabla }%
:T_{r}^{s}\rightarrow T^{\ast }\otimes T_{r}^{s}.$

We now have the following important

\begin{proposition}
The following are equivalent

$i)$ $\xi \in T_{r}^{s}$ is $\varepsilon $-invariant

$ii)$ $\xi $ is $\widetilde{\cdot }$ -parallel, that is, $\widetilde{\nabla }%
_{X}\xi =0$ for all $X\in \mathfrak{X}(M)$

If $X\in \widehat{\Theta }$, then $\mathcal{L}_{X}=\widetilde{\nabla }_{X}.$
In particular, if $\widehat{\mathfrak{R}}=0$, then $i)$ and $ii)$ are
equivalent to $iii)$ $\mathcal{L}_{X}\xi =0$ for all $X\in \widehat{\Theta }$
\end{proposition}

$i)\Rightarrow ii):$ When we differentiate (8) at $x=p$ and substitute from
(2), we get $(\widetilde{\nabla }\xi )(p)=0$ and therefore $\widetilde{%
\nabla }\xi =0$ since $p$ is arbitrary.

$ii)\Rightarrow i):$ If $\widetilde{\nabla }\xi =0,$ we claim that $\xi $ is
determined near $p$ by $\xi (p)$ (note that $\widetilde{\nabla }\xi =0$ may
not have any solution other than $\xi =0$ because we are not making any
assumption on $\widehat{\mathfrak{R}}$ and/or $T).$ Indeed the restriction
of a solution $\xi $ to a smooth path starting from $p$ is a solution of a
linear $ODE$ with initial condition $\xi (p).$ Therefore the solution is
unique on this path and hence near $p.$ Therefore the solution $\xi $ must
coincide with the $\varepsilon $-invariant extension of $\xi (p)$ since they
both solve $\widetilde{\nabla }\xi =0$ with the same initial condition.

$ii)\Leftrightarrow iii):$ Let $X\in \widehat{\Theta },$ $\eta \in \mathfrak{%
X}(M).$ We have

\begin{eqnarray}
\left( \mathcal{L}_{X}\eta \right) ^{i} &=&[X,\eta ]^{i}=X^{a}\frac{\partial
\eta ^{i}}{\partial x^{a}}-\eta ^{a}\frac{\partial X^{i}}{\partial x^{a}} \\
&=&X^{a}\frac{\partial \eta ^{i}}{\partial x^{a}}-\eta ^{a}\Gamma
_{ba}^{i}X^{b}  \notag \\
&=&\left( \frac{\partial \eta ^{i}}{\partial x^{b}}-\eta ^{a}\Gamma
_{ba}^{i}\right) X^{b}  \notag \\
&=&\widetilde{\nabla }_{X}(\eta )^{i}  \notag
\end{eqnarray}

Thus the operators $\widetilde{\nabla }_{X}$ and $\mathcal{L}_{X}$ coincide
on vector fields and we conclude that they coincide on all tensor fields.

The elementary computation in (10) is extremely important for it shows that $%
\widetilde{\nabla }_{X}$ is actually a substitute for the real object $%
\mathcal{L}_{X}$ which emerges fully when $\widehat{\mathfrak{R}}=0.$

It is useful to recast some of the above local formulas in a coordinate free
language.

As already noted above, we have the operator $\widetilde{\nabla }%
_{X}:T\rightarrow T$ defined by $\left( \widetilde{\nabla }_{X}Y\right) ^{i}%
\overset{def}{=}\left( \frac{\partial Y^{i}}{\partial x^{a}}-\Gamma
_{ab}^{i}Y^{b}\right) X^{a}.$ Therefore $Y\in \widetilde{\Theta }%
\Leftrightarrow \widetilde{\nabla }_{X}Y=0$ for all $X\in \mathfrak{X}(M).$
Similarly we define $\widehat{\nabla }_{X}$ and deduce the similar
statements. We easily deduce the formulas

\begin{eqnarray}
\widetilde{\nabla }_{X}Y-\widehat{\nabla }_{X}Y &=&T(X,Y) \\
\widetilde{\nabla }_{X}Y-\widetilde{\nabla }_{Y}X &=&[X,Y]+T(X,Y)  \notag \\
\widehat{\nabla }_{X}Y-\widetilde{\nabla }_{Y}X &=&[X,Y]  \notag \\
\widetilde{\nabla }_{X}Y+\widetilde{\nabla }_{Y}X &=&\widehat{\nabla }_{X}Y+%
\widehat{\nabla }_{Y}X  \notag
\end{eqnarray}

If $Y\in \widetilde{\Theta }$ and $X\in \widehat{\Theta }$, then $%
T(X,Y)=[X,Y]=0.$ Also, $\widetilde{\mathfrak{R}}\in \Lambda ^{2}(T^{\ast
})\otimes T^{\ast }\otimes T$ defined by $\left( \widetilde{\mathfrak{R}}%
(X,Y)(Z)\right) ^{i}\overset{def}{=}\widetilde{\mathfrak{R}}%
_{ab,c}^{i}X^{a}Y^{b}Z^{c}$ and $\widehat{\mathfrak{R}}$ is similarly
defined.

Now let $(M,\varepsilon )$ be a parallelizable manifold with the splitting $%
\varepsilon .$ If $X,Y\in \widetilde{\Theta }$, by direct computation we
check the identity

\begin{equation}
\widetilde{\nabla }_{Z}[X,Y]=\widehat{\mathfrak{R}}(X,Y)Z\text{ \ \ }Z\in
\mathfrak{X}(M)
\end{equation}

Therefore $[X,Y]\in \widetilde{\Theta }$ for all $X,Y\in \widetilde{\Theta }%
\Leftrightarrow \widetilde{\nabla }_{Z}[X,Y]=0$ for all $Z\in \mathfrak{X}%
(M)\Leftrightarrow \widehat{\mathfrak{R}}(X,Y)=0$ for all $X,Y\in \widetilde{%
\Theta }\Leftrightarrow \widehat{\mathfrak{R}}(p)(X_{p},Y_{p})=0$ for all $%
p\in M$ and $X_{p},Y_{p}\in T_{p}=$ the tangent space at $p\Leftrightarrow
\widehat{\mathfrak{R}}=0.$

Now we make the following important

\begin{definition}
For $X,Y\in \mathfrak{X}(M),$ we call $T(X,Y)\in \mathfrak{X}(M)$ the
algebraic bracket of $X,Y.$
\end{definition}

By (6), we have

\begin{equation}
\left( T(X,Y)(p)\right) ^{i}=T(p)_{ab}^{i}X^{a}Y^{b}
\end{equation}

It is important to observe that the algebraic bracket $T(X,Y)$ is defined on
any parallelizable manifold for we need only (2) to define it. If $\widehat{%
\mathfrak{R}}=0,$ then $\widetilde{\Theta }$ is a Lie algebra and (11) shows
$T(X,Y)=-[X,Y]$ for $X,Y\in \widetilde{\Theta }$. So $\widehat{\mathfrak{R}}%
=0$ forces $T(\cdot ,\cdot )$ to satisfy the Jacobi identity which need not
hold in general$.$ For simplicity of notation, we set

\begin{equation}
J(X,Y,Z)\overset{def}{=}T(X,T(Y,Z))+T(Z,T(X,Y))+T(Y,T(Z,X))
\end{equation}%
so that $J$ is an $3$-form and the Jacobi identity is equivalent to $J=0.$

The following proposition is as fundamental as the $LTFT.$

\begin{proposition}
On a parallelizable manifold $(M,\varepsilon )$, the following identities
hold.
\end{proposition}

\begin{equation}
(\widetilde{\nabla }_{X}T)(Y,Z)+(\widetilde{\nabla }_{Z}T)(X,Y)+(\widetilde{%
\nabla }_{Y}T)(Z,X)=J(X,Y,Z)
\end{equation}

\begin{equation}
\widehat{\mathfrak{R}}(X,Y)(Z)=\text{ }(\widetilde{\nabla }_{X}T)(Z,Y)+(%
\widetilde{\nabla }_{Y}T)(X,Z)+J(X,Y,Z)\text{\ }
\end{equation}%
where $X,Y,Z\in \mathfrak{X}(M).$

In coordinates, (15) reduces to checking the identity

\begin{equation}
\left[ \frac{\partial T_{kj}^{i}}{\partial x^{r}}-\Gamma
_{ra}^{i}T_{kj}^{a}+\Gamma _{rk}^{a}T_{aj}^{i}+\Gamma _{rj}^{a}T_{ka}^{i}%
\right]
_{[rkj]}=T_{rk}^{a}T_{aj}^{i}+T_{jr}^{a}T_{ak}^{i}+T_{kj}^{a}T_{ar}^{i}
\end{equation}

To prove (16) we proceed as

\begin{eqnarray}
\widehat{\mathfrak{R}}_{rj,k}^{i} &=&\widehat{\mathfrak{R}}_{rj,k}^{i}-%
\widetilde{\mathfrak{R}}_{rj,k}^{i}  \notag \\
&=&\left[ \frac{\partial T_{kj}^{i}}{\partial x^{r}}-\Gamma
_{ra}^{i}T_{kj}^{a}+\Gamma _{rk}^{a}T_{aj}^{i}+\Gamma _{rj}^{a}T_{ka}^{i}%
\right] _{[rj]} \\
&&+T_{aj}^{i}T_{rk}^{a}+T_{ar}^{i}T_{kj}^{a}+T_{ak}^{i}T_{jr}^{a}  \notag
\end{eqnarray}%
and check the second equality in (18) using (17).

\begin{corollary}
$\widehat{\mathfrak{R}}(X,Y)(Z)=(\widetilde{\nabla }_{Z}T)(X,Y)$ (which
generalizes (12))$.$ In particular, $(M,\varepsilon )$ is a local Lie group
if and only if $\widetilde{\nabla }T=0.$
\end{corollary}

Note that $(\widetilde{\nabla }_{Z}T)(X,Y)$ and $\widetilde{\nabla }%
_{Z}\left( T(X,Y)\right) $ are different and related by $(\widetilde{\nabla }%
_{Z}T)$ $(X,Y)=\widetilde{\nabla }_{Z}\left( T(X,Y)\right) +T(\widetilde{%
\nabla }_{Z}X,Y)+T(X,\widetilde{\nabla }_{Z}Y).$

The above corollary shows that curvature is determined by torsion. Therefore
the hero of the present tensor calculus is torsion and not curvature. This
will be further supported by the \textit{even }dimensional characteristic
classes constructed at the end of Section 4.

\begin{corollary}
(Bianchi identity) $\widehat{\mathfrak{R}}(X,Y)(Z)+\widehat{\mathfrak{R}}%
(Z,X)(Y)+\widehat{\mathfrak{R}}(Y,Z)(X)=J(X,Y,Z)$
\end{corollary}

Now we define the operator $[\widetilde{\nabla }]:\Lambda ^{k}(T^{\ast
})\otimes T\rightarrow \Lambda ^{k+1}(T^{\ast })\otimes T$ by

\begin{eqnarray}
&&\left( \lbrack \widetilde{\nabla }]\xi \right) _{rst....m}^{i}\overset{def}%
{=}\left[ \widetilde{\nabla }_{r}\xi _{st...m}^{i}\right] _{[rst...m]} \\
&=&\widetilde{\nabla }_{r}\xi _{st...m}^{i}-\widetilde{\nabla }_{s}\xi
_{rt...m}^{i}-\widetilde{\nabla }_{t}\xi _{sr...m}^{i}-...\widetilde{\nabla }%
_{m}\xi _{st...r}^{i}  \notag
\end{eqnarray}%
or equivalently

\begin{eqnarray*}
([\widetilde{\nabla }]\xi )(X_{0},Y_{1},...,X_{k}) &=&(\widetilde{\nabla }%
_{X_{0}}\xi )(X_{1},X_{2},...,X_{k})-(\widetilde{\nabla }_{X_{1}}\xi
)(X_{0},X_{2},...,X_{k}) \\
&&-(\widetilde{\nabla }_{X_{2}}\xi )(X_{1},X_{0},...,X_{k})-...(\widetilde{%
\nabla }_{X_{0}}\xi )(X_{1},X_{2},...,X_{k}
\end{eqnarray*}%
and obtain the sequence

\begin{equation}
0\rightarrow \widetilde{\Theta }\rightarrow T\overset{[\widetilde{\nabla }]}{%
\longrightarrow }T^{\ast }\otimes T\overset{[\widetilde{\nabla }]}{%
\longrightarrow }\Lambda ^{2}(T^{\ast })\otimes T\overset{[\widetilde{\nabla
}]}{\longrightarrow }....\overset{[\widetilde{\nabla }]}{\longrightarrow }%
\Lambda ^{n}(T^{\ast })\otimes T
\end{equation}

Unfortunately (20) is not a complex even though $\widetilde{\mathfrak{R}}=0$
for $T$ is the obstruction to its exactness. For instance, we have $%
\widetilde{\nabla }_{k}\widetilde{\nabla }_{j}\xi ^{i}-\widetilde{\nabla }%
_{j}\widetilde{\nabla }_{k}\xi ^{i}=T_{kj}^{a}\widetilde{\nabla }_{a}\xi
^{i} $ (this formula is written incorrectly in [1]). But (6) shows that $%
T=0\Leftrightarrow $ $\widetilde{\Theta }\simeq \widetilde{\Theta }_{p}$ is
abelian and therefore $T=0$ is a very strong assumption. However, as we will
see in the next section, a fine resolution of the sheaf $\widetilde{\Theta }$
is an object quite different from (20) but related to (20) in a subtle way.

\section{Some isomorphic cohomology groups}

In geometry it is standard to define some cohomology groups in certain
situations. These are shortly

1) For a sheaf $\Theta $ on $M$, we have the sheaf cohomology groups.

2) For a vector bundle $E\rightarrow M$ with a flat connection, we have the
complex of the $E$-valued forms on $M.$

3) For a vector bundle $E\rightarrow M$ and an involutive sytem $R$ on $%
J_{k}E\rightarrow M$, we have the cohomology of the linear Janet sequence of
$R$ ([11]).

4) For a Lie algebra $\mathfrak{g}$\ with a representation on a vector space
$V$, we have the Chevalley-Eilenberg cohomology groups ([2])

5) For an algebroid $A\rightarrow M$ with a representation on a vector
bundle $T\rightarrow M,$ we have the complex of $T$-valued $A$-forms ([3],
[8])

In this section, we will show that the first four of the above general
constructions coincide in a special case and we will comment on the last one
at the end of this section. Probably 3) is the least known among the others
so we will start with 3) to advertise it.

Let $(M,\varepsilon )$ be a local Lie group. Since $\widehat{\mathfrak{R}}%
=0, $ the linear system $R\overset{def}{=}$\ $\widehat{\Gamma }(T)\subset T$
is involutive where $T\rightarrow M$ is the tangent bundle. We have the
operator $\widehat{\nabla }:T\rightarrow T^{\ast }\otimes T$ and $\widehat{%
\Theta }=Ker(\widehat{\nabla })=$ the sheaf of local solutions of $2^{\prime
}).$ Thus we have the exact sequence $0\rightarrow \widehat{\Theta }%
\rightarrow T\overset{\widehat{\nabla }}{\rightarrow }T^{\ast }\otimes T.$
To construct the Janet sequence of $R=$\ $\widehat{\Gamma }(T),$ we check
the integrability conditions of $\widehat{\nabla }.$ So we set

\begin{equation}
\frac{\partial \xi ^{i}}{\partial x^{j}}-\Gamma _{aj}^{i}\xi ^{a}=\eta
_{j}^{i}
\end{equation}

We differentiate (21) with respect to $x^{k},$ alternate $k,j$ and
substitute $\frac{\partial \xi ^{a}}{\partial x^{k}}$ back from (21).
Arranging terms, we get
\begin{equation}
0=\widehat{\mathcal{R}}_{jk,a}^{i}\xi ^{a}=\frac{\partial \eta _{j}^{i}}{%
\partial x^{k}}-\frac{\partial \eta _{k}^{i}}{\partial x^{j}}+\Gamma
_{aj}^{i}\eta _{k}^{a}-\Gamma _{ak}^{i}\eta _{j}^{a}
\end{equation}

Denoting $\widehat{\nabla }$ by $\widehat{d}_{0},$ we now define $\widehat{%
d_{1}}:T^{\ast }\otimes T\rightarrow \Lambda ^{2}(T^{\ast })\otimes T$ by
the formula (22), that is, $\widehat{d}_{1}:(\eta _{j}^{i})\rightarrow \left[
\widehat{d}_{k}\eta _{j}^{i}\right] _{[k,j]}$ where $d_{k}\eta _{j}^{i}%
\overset{def}{=}\frac{\partial \eta _{j}^{i}}{\partial x^{k}}-\Gamma
_{ak}^{i}\eta _{j}^{a}.$ By definition $\widehat{d}_{1}\circ \widehat{d}%
_{0}=0.$ We now check the integrability conditions of the operator $\widehat{%
d}_{1}$. So we set

\begin{equation}
\frac{\partial \eta _{j}^{i}}{\partial x^{k}}-\frac{\partial \eta _{k}^{i}}{%
\partial x^{j}}+\Gamma _{aj}^{i}\eta _{k}^{a}-\Gamma _{ak}^{i}\eta
_{j}^{a}=\phi _{kj}^{i}
\end{equation}

We now differentiate (23) with respect to $x^{r}$ and alternate $r$ with $%
k,j.$ We substitute three differences each of the form $\frac{\partial \eta
_{k}^{a}}{\partial x^{r}}-\frac{\partial \eta _{r}^{a}}{\partial x^{k}}$
back from (23) and transpose all the terms which contain the components of $%
\phi $ to the right of (23) where we get $\left[ \widehat{d}_{r}\phi
_{kj}^{i}\right] _{[rkj]}=\widehat{d}_{r}\phi _{kj}^{i}-\widehat{d}_{k}\phi
_{rj}^{i}-\widehat{d}_{j}\phi _{kr}^{i}$ with $\widehat{d}_{r}\phi _{kj}^{i}%
\overset{def}{=}\frac{\partial \phi _{kj}^{i}}{\partial x^{r}}-\Gamma
_{ar}^{i}\phi _{kj}^{a}.$ On the left we are left with three terms each of
the form $\widehat{\mathfrak{R}}_{rk,a}^{i}\eta _{j}^{a}$ so this sum
vanishes since $\widehat{\mathfrak{R}}=0$. Thus we get the operator $%
\widehat{d}_{2}:\Lambda ^{2}(T^{\ast })\otimes T\rightarrow \Lambda
^{3}(T^{\ast })\otimes T$ defined by $(\phi _{kj}^{i})\rightarrow \left[
\widehat{d}_{r}\phi _{kj}^{i}\right] _{[rkj]}$ and again by definition $%
\widehat{d}_{2}\circ \widehat{d}_{1}=0$. Iterating this process of checking
"the integrability conditions of the integrability conditions" gives the
linear Janet sequence of the system $R$ which is

\begin{equation}
0\rightarrow \widehat{\Theta }\rightarrow T\overset{\widehat{d}}{%
\longrightarrow }T^{\ast }\otimes T\overset{\widehat{d}}{\longrightarrow }%
\Lambda ^{2}(T^{\ast })\otimes T\overset{\widehat{d}}{\longrightarrow }....%
\overset{\widehat{d}}{\longrightarrow }\Lambda ^{n}(T^{\ast })\otimes T
\end{equation}%
and the operator $\widehat{d}_{k}:\Lambda ^{k}(T^{\ast })\otimes
T\rightarrow \Lambda ^{k+1}(T^{\ast })\otimes T$ is defined by

\begin{equation}
\left[ \widehat{d}_{r}\xi _{jl....m}^{i}\right] _{[rjl...m]}\overset{def}{=}%
\widehat{d}_{r}\xi _{jl....m}^{i}-\widehat{d}_{j}\xi _{rl....lm}^{i}-%
\widehat{d}_{l}\xi _{jr....m}^{i}-....-\widehat{d}_{m}\xi _{jl....r}^{i}
\end{equation}%
where $\widehat{d}_{r}\omega _{jl....m}^{i}\overset{def}{=}\frac{\partial
\xi _{jl....m}^{i}}{\partial x^{r}}-\Gamma _{ar}^{i}\xi _{jl....m}^{a}.$

The Janet sequence (24) is locally exact. Since partition of unity applies
to the sections of the spaces in (24), (24) is a fine resolution of the
sheaf $\widehat{\Theta }.$ Further, $2^{\prime })$ defines a connection on
the tangent bundle $T\rightarrow M$ which is flat since $\widehat{\mathfrak{R%
}}=0$ and our local formulas define the "exterior covariant
differentiation". Using a suggestive notation, we therefore have \

\begin{proposition}
$H^{\ast }(\widehat{\Theta })_{sheaf}\simeq H^{\ast }(\widehat{\Theta }%
)_{Janet}\simeq H^{\ast }(\widehat{\Theta })_{flat}$
\end{proposition}

Note that we did not make use of the Lie algebra structure of $\widehat{%
\Theta }$ yet. Thus we could start with $2)$ instead of $2^{\prime })$ and
construct (24) for the sheaf $\widetilde{\Theta }$ with the operators $%
\widetilde{d}$ accordingly defined and still get Proposition 6 with $%
\widetilde{\Theta }$ replacing $\widehat{\Theta }$. This latter construction
works for any parallelizable $(M,\varepsilon )$ since we always have $%
\widetilde{\mathfrak{R}}=0.$ Note that the operators $\widetilde{d},$ $[%
\widetilde{\nabla }]:\Lambda ^{k}(T^{\ast })\otimes T\rightarrow \Lambda
^{k+1}(T^{\ast })\otimes T$ coincide only for $k=1.$

Now we compare (19) and (25). We have
\begin{eqnarray}
\widehat{d}_{r}\xi _{jl....m}^{i}+T_{ar}^{i}\xi _{jl....m}^{a} &=&\frac{%
\partial \xi _{jl....m}^{i}}{\partial x^{r}}-\Gamma _{ar}^{i}\xi
_{jl....m}^{a}+T_{ar}^{i}\xi _{jl....m}^{a}  \notag \\
&=&\frac{\partial \xi _{jl....m}^{i}}{\partial x^{r}}-\Gamma _{ra}^{i}\xi
_{jl....m}^{a}
\end{eqnarray}

So if we add $\left[ T_{ar}^{i}\omega _{jl....m}^{a}\right] _{[rjl...m]}$ to
(25), this will change the term $\Gamma _{ar}^{i}\xi _{jl....m}^{a}$ in $%
\widehat{d}$ to $\Gamma _{ra}^{i}\xi _{jl....m}^{a}$ in $[\widetilde{\nabla }%
]$ but we must add also the remaining term

\begin{equation*}
\left[ \Gamma _{rj}^{a}\omega _{al...m}^{i}+\Gamma _{rl}^{a}\omega
_{ja...m}^{i}...+\Gamma _{rm}^{a}\omega _{jl...a}^{i}\right] _{[rjl...m]}
\end{equation*}%
to $\widehat{d}$ to get $[\widetilde{\nabla }].$ So we deduce

\begin{eqnarray}
([\widetilde{\nabla }]\xi )_{rjl...m}^{i} &=&(\widehat{d}\xi )_{rjl...m}^{i}+%
\left[ T_{ar}^{i}\xi _{jl...m}^{a}\right] _{[rjl...m]} \\
&&+\left[ \Gamma _{rj}^{a}\xi _{al...m}^{i}+\Gamma _{rl}^{a}\xi
_{ja...m}^{i}...+\Gamma _{rm}^{a}\xi _{jl...a}^{i}\right] _{[rjl...m]}
\notag
\end{eqnarray}%
\qquad

We call some $\xi \in \Lambda ^{k}(T^{\ast })\otimes T$ shortly a $k$-form.
Now let $\widetilde{\Lambda ^{k}(T^{\ast })\otimes T}$ denote the space of $%
\varepsilon $-invariant (or $\widetilde{\cdot }$ -parallel in view of
Proposition 1) $k$-forms so that $\widetilde{\Theta }=\widetilde{T}$.

We now have the following important

\begin{proposition}
$\widehat{d}:\widetilde{\Lambda ^{k}(T^{\ast })\otimes T}$ $\longrightarrow $
$\widetilde{\Lambda ^{k+1}(T^{\ast })\otimes T}$
\end{proposition}

Indeed, if $\xi \in \widetilde{\Lambda ^{k}(T^{\ast })\otimes T},$ then $[%
\widetilde{\nabla }]\xi =0$ by Proposition 1 and (27) gives

\begin{equation}
(\widehat{d}\xi )_{rjl...m}^{i}=-\left[ T_{ar}^{i}\xi _{jl...m}^{a}\right]
_{[rjl...m]}-\left[ \Gamma _{rj}^{a}\xi _{al...m}^{i}+\Gamma _{rl}^{a}\xi
_{ja...m}^{i}...+\Gamma _{rm}^{a}\xi _{jl...a}^{i}\right] _{[rjl...m]}
\end{equation}

By Corollary 4 $\widetilde{\nabla }T=0$. Since products and contractions of $%
\widetilde{\cdot }$ -parallel tensors are $\widetilde{\cdot }$ -parallel,
the first term on the RHS of (28) is $\widetilde{\cdot }$ -parallel. Note
that the expression inside the paranthesis of the second term is already
alternating in $jl...m.$ It is not difficult to show that when we alternate
also $r,$ the second term becomes a sum of terms of the form $T_{\cdot \cdot
}^{a}\xi _{a\cdot \cdot ...\cdot }$ and is therefore also $\widetilde{\cdot }
$ -parallel finishing the proof. For later use, we list the first three of
these alternations:

\begin{eqnarray}
\left[ \Gamma _{rk}^{a}\xi _{a}^{i}\right] _{[rk]} &=&T_{rk}^{a}\xi _{a}^{i}
\notag \\
\left[ \Gamma _{rk}^{a}\xi _{aj}^{i}+\Gamma _{rj}^{a}\xi _{ka}^{i}\right]
_{[rkj]} &=&T_{rk}^{a}\xi _{aj}^{i}+T_{jr}^{a}\xi _{ak}^{i}+T_{kj}^{a}\xi
_{ar}^{i} \\
\left[ \Gamma _{rk}^{a}\xi _{ajm}+\Gamma _{rj}^{a}\xi _{kam}+\Gamma
_{rm}^{a}\xi _{kja}\right] _{[rkjm]} &=&T_{rk}^{a}\xi _{ajm}+T_{jr}^{a}\xi
_{akm}+T_{rm}^{a}\xi _{akj}  \notag \\
&&+T_{mj}^{a}\xi _{akr}+T_{km}^{a}\xi _{ajr}+T_{kj}^{a}\xi _{arm}  \notag
\end{eqnarray}

It is worthwhile to write (27) in a coordinate free form but this would
carry us away from our main purpose here. Proposition 7 gives the subcomplex

\begin{equation}
\widetilde{\Theta }=\widetilde{T}\overset{\widehat{d}}{\longrightarrow }%
\widetilde{T^{\ast }\otimes T}\overset{\widehat{d}}{\longrightarrow }%
\widetilde{\Lambda ^{2}(T^{\ast })\otimes T}\overset{\widehat{d}}{%
\longrightarrow }....\overset{\widehat{d}}{\longrightarrow }\widetilde{%
\Lambda ^{n}(T^{\ast })\otimes T}
\end{equation}%
We denote the cohomology of (30) by $H^{\ast }(\widetilde{\Theta },%
\widetilde{\Theta })$ for the reason which will become clear shortly. The
inclusion of (30) into (24) does not induce isomorphism in cohomology in
general. Recal that the cohomology of invariant scalar valued forms on a Lie
group is isomorphic to the de Rham cohomology if $M$ is compact. The well
known "averaging with respect to the Haar measure" proof works also for
forms with values in a vector bundle if the the group permutes the fibers of
this vector bundle in a way consistent with its action on the base manifold
as in the case of (24). Thus we obtain the following fundamental

\begin{proposition}
If $M$ is compact, then $H^{\ast }(\widetilde{\Theta },\widetilde{\Theta })$
is isomorphic to the cohomology groups in Proposition 6 in positive degrees.
\end{proposition}

Note that Propositions 7, 8 imply $H^{\ast }(\widetilde{\Theta },\widetilde{%
\Theta })=H^{\ast }(\widehat{\Theta })$ for compact $M$ but Proposition 9
below will show that this confusion with notation does not stem from us.

Now we fix some base point $p\in M$ and consider the vector space $\Lambda
^{k}(T_{p}^{\ast })\otimes T_{p}$ where $T_{p}$ and $T_{p}^{\ast }$ are the
tangent and cotangent spaces at $p.$ We define a map $\Lambda
^{k}(T_{p}^{\ast })\otimes T_{p}\rightarrow \Lambda ^{k+1}(T_{p}^{\ast
})\otimes T_{p}$ as follows: if $\sigma _{p}\in \Lambda ^{k}(T_{p}^{\ast
})\otimes T_{p}$, we extend $\sigma _{p}$ to an element $\sigma $ of $%
\widetilde{\Lambda ^{k}(T^{\ast })\otimes T},$ apply $\widehat{d}$ to $%
\sigma $ and evaluate the result at $p.$ We denote this operator by $%
\widehat{d}(p)$ and call it the localization of $\widehat{d}.$ So we have $%
\widehat{d}(p)(\sigma _{p})\overset{def}{=}$ $(\widehat{d}\sigma )(p)$. In
this way we obtain the purely algebraic complex

\begin{equation}
T_{p}\overset{\widehat{d}(p)}{\longrightarrow }T_{p}^{\ast }\otimes T_{p}%
\overset{\widehat{d}(p)}{\longrightarrow }\Lambda ^{2}(T_{p}^{\ast })\otimes
T_{p}\overset{\widehat{d}(p)}{\longrightarrow }....\overset{\widehat{d}(p)}{%
\longrightarrow }\Lambda ^{n}(T_{p}^{\ast })\otimes T_{p}
\end{equation}

We denote the cohomology of (31) by $H^{\ast }(\widetilde{\Theta }_{p},%
\widetilde{\Theta }_{p}).$ Since the evaluation and its inverse extension
are Lie algebra isomorphisms, $H^{\ast }(\widetilde{\Theta },\widetilde{%
\Theta })\simeq H^{\ast }(\widetilde{\Theta }_{p},\widetilde{\Theta }_{p}).$
In particular $H^{\ast }(\widetilde{\Theta }_{p},\widetilde{\Theta }%
_{p})\simeq H^{\ast }(\widetilde{\Theta }_{q},\widetilde{\Theta }_{q})$ and $%
\simeq $ is canonical: since $\widehat{\mathfrak{R}}=0\Rightarrow \mathcal{%
R(\varepsilon )=}0,$ there is a unique local diffeomorphism $f$ of the
pseudogroup $\mathcal{G}$ with $f(p)=q.$ Since $f$ commutes with the
operators in (24), it induces $\simeq .$

Now our purpose is to explicitly compute $H^{i}(\widetilde{\Theta }_{p},%
\widetilde{\Theta }_{p})$ for $i=0,1.$ So let $\xi _{p},\eta _{p}\in
\widetilde{\Theta }_{p}.$ By definition $\widehat{d}(p)(\xi _{p})=(\widehat{d%
}\xi )(p)$ where $\xi $ is the $\varepsilon $-invariant extension of $\xi
_{p}$ as above. Now $(\widehat{d}\xi )(\eta )=[\xi ,\eta ]$ by (10).
Therefore $\widehat{d}(p)(\xi _{p})=0\Leftrightarrow \lbrack \xi ,\eta
](p)=[\xi _{p},\eta _{p}]_{p}=0$ for all $\eta _{p}\in \widetilde{\Theta }%
_{p}\Leftrightarrow \xi _{p}\in C(\widetilde{\Theta }_{p})=$ the center of $%
\widetilde{\Theta }_{p}.$ Therefore $H^{0}(\widetilde{\Theta }_{p},%
\widetilde{\Theta }_{p})\simeq C(\widetilde{\Theta }_{p}).$

Now suppose $\widehat{d}(p)(\xi _{p})=0$ for some $\xi _{p}=((\xi
_{p})_{j}^{i})\in T_{p}^{\ast }\otimes T_{p}.$ The formula (28) specializes
to

\begin{eqnarray}
(\widehat{d}\xi )_{rj}^{i} &=&-\left[ T_{ra}^{i}\xi _{j}^{a}\right] _{[rj]}-%
\left[ \Gamma _{rj}^{a}\xi _{a}^{i}\right] _{[rj]}  \notag \\
&=&-T_{rj}^{a}\xi _{a}^{i}-T_{ra}^{i}\xi _{j}^{a}+T_{ja}^{i}\xi _{r}^{a}
\end{eqnarray}

Recall that a linear map $\xi _{p}:\widetilde{\Theta }_{p}\rightarrow
\widetilde{\Theta }_{p}$ is a derivation if

\begin{equation}
\xi _{p}[\sigma _{p},\eta _{p}]_{p}-[\xi _{p}\sigma _{p},\eta
_{p}]_{p}-[\sigma _{p},\xi _{p}\eta _{p}]_{p}=0
\end{equation}%
for all $\sigma _{p},\eta _{p}$. Now a straightforward computation using (6)
shows that with $\sigma _{p}=\left[ \frac{\partial }{\partial x^{r}}\right]
_{p}$ and $\eta _{p}=\left[ \frac{\partial }{\partial x^{j}}\right] _{p},$
the LHS of (33) coincides with the evaluation of the RHS of (32) at $p.$
Therefore $1$-cocycles in (31) are precisely derivations and (10) shows that
a derivation is a boundary if and only if it is inner. Therefore $H^{1}(%
\widetilde{\Theta }_{p},\widetilde{\Theta }_{p})\simeq Der(\widetilde{\Theta
}_{p})/Inn(\widetilde{\Theta }_{p}).$

Not surprisingly, we now have

\begin{proposition}
(31) is the complex which computes the cohomology of the Lie algebra $%
\widetilde{\Theta }_{p}$ with respect to its adjoint representation on
itself, that is, the deformation cohomology $H^{\ast }(\widetilde{\Theta }%
_{p},\widetilde{\Theta }_{p}).$
\end{proposition}

We will omit the straightforward verification.

Now the algebroid $T\rightarrow M$ has a representation on the vector bundle
$T\rightarrow M$ defined by the ordinary bracket of vector fields. A $T$
valued $k$-form is an element of $\Lambda ^{k}(T^{\ast })\otimes T$ and we
get the same spaces in (24). Let us denote the cohomology of this complex by
$H^{\ast }(T,T)$. The first operator $\delta :T\rightarrow T^{\ast }\otimes
T $ in this complex is defined by

\begin{equation}
(\delta )(\eta )(\xi )=[\eta ,\xi ]=\mathcal{L}_{\eta }\xi \text{ \ \ }\eta ,%
\text{ }\xi \in \mathfrak{X}(M)
\end{equation}%
and it is worthwhile to observe the difference between (10) and (34). It is
easy to show that $Ker(\delta )=\{0\}.$ Therefore, the complex computing $%
H^{\ast }(T,T)$ is either locally exact in which case $H^{\ast }(T,T)$
vanishes since this complex is a fine resolution of the sheaf $\{0\}$ or it
is not locally exact. In view of [1], Section 5, we believe that the second
is the case.

In [1] we defined some odd degree characteristic classes in the Lie algebra
cohomology and therefore also in the de Rham cohomology of a local Lie group
$(M,\varepsilon ).$ In essence these classes are obtained by interpreting
torsion as a $1$-form, taking its powers and then taking its trace (see [1]
for details). The next proposition shows that the torsion itself as a $2$%
-form defines a closed form in the comlex (30).

\begin{proposition}
$\widehat{d}T=0$
\end{proposition}

The reason is that $T$ is $\widetilde{\cdot }$ -parallel by Corollary 4 and
the LHS of (27) vanishes so that we have (28) with $k=2$ and $\xi =T$. Now
both paranthesis on the RHS of (28) boil down to $J(X,Y,Z).$ Thus the RHS of
(28) is $-2J(X,Y,Z)$ which vanishes by (15).

So we obtain the characteristic class $T^{1}\in H^{2}(\widetilde{\Theta },%
\widetilde{\Theta }).$ Multiplying $T$'s $k$-times and alternating, we get
the $(k+1)$ form $T^{k}.$ For instance

\begin{eqnarray}
&&(T^{2})_{jkl}^{i}\overset{def}{=}\left[ T_{ja}^{i}T_{kl}^{a}\right]
_{[jkl]} \\
&&(T^{3})_{jklm}^{i}\overset{def}{=}\left[ T_{ja}^{i}T_{kb}^{a}T_{lm}^{b}%
\right] _{[jklm]}  \notag
\end{eqnarray}
\

Note that $T^{2}(X,Y,Z)=J(X,Y,Z)=0.$ Similarly, all the odd degree forms in
(35) vanish and using Proposition 10, it is not difficult to show that even
degree forms are closed in (30). Thus we obtain the cohomology classes

\begin{equation}
\lbrack T^{2k-1}]\in H^{2k}(\widetilde{\Theta },\widetilde{\Theta })
\end{equation}

Now $Trace$ defines a map $tr$ from the complex (30) to the complex
computing the Lie algebra cohomology of $\widetilde{\Theta }$ but with one
shift in degrees, that is

\begin{equation}
\begin{array}{ccccccccc}
\widetilde{\Theta }=\widetilde{T} & \overset{\widehat{d}}{\longrightarrow }
& \widetilde{T^{\ast }\otimes T} & \overset{\widehat{d}}{\longrightarrow } &
\widetilde{\Lambda ^{2}(T^{\ast })\otimes T} & \overset{\widehat{d}}{%
\longrightarrow } & \widetilde{\Lambda ^{3}(T^{\ast })\otimes T} & \overset{%
\widehat{d}}{\longrightarrow } & .... \\
&  & \downarrow tr &  & \downarrow tr &  & \downarrow tr &  &  \\
&  & \mathbb{R} & \overset{d}{\longrightarrow } & \widetilde{\Lambda
^{1}(T^{\ast })} & \overset{d}{\longrightarrow } & \widetilde{\Lambda
^{2}(T^{\ast })} & \overset{d}{\longrightarrow } & ....%
\end{array}%
\end{equation}%
For instance, $tr(T^{1})=T_{ja}^{a}$ which is the $1$-form in [1] and $%
tr(T^{3})=(T^{3})_{jkla}^{a}$ which is the $3$-form in [1] and we easily see
that $tr$ maps the even degree characteristic classes (35) to the odd degree
secondary characteristic classes in [1]. However, note that $[T^{2k-1}]$ are
also secondary for they are defined only on a local Lie group, that is, when
$\widehat{\mathfrak{R}}=0\mathfrak{.}$ We do not know the geometric meaning
of these classes.

\section{Gauge transformations}

This paper was planned to end as a short note at the end of Section 4 but
Proposition 9 forced us to go further. Indeed, the study of deformations of
algebraic structures is initiated in [4] shortly after [7] and later studied
by various authors. However, for someone not familiar with this algebraic
theory, at first sight it is not clear at all why (31) is called
"deformation complex". What are we really deforming? As explained in Section
2, on a local Lie group $(M,\varepsilon )$ the vector fields in $\widehat{%
\Theta }$ integrate to the pseudogroup $\mathcal{G}$ whose $1$-arrows are
given by $\varepsilon $ that we start with. Further if $(M,\varepsilon )$ is
complete and $M$ is simply connected, $\mathcal{G}$ globalizes to a
transformation group which acts simply transitively on $M$ with
infinitesimal generators $\widehat{\Theta }.$ Clearly, the deformation of $%
\mathcal{G}$ with time (to be defined in this section) is a geometrically
very intuitive concept and naturally the problem arises how to deduce
algebraic deformation from this geometric deformation. This deduction will
be our main concern in the rest of this note.

We recall the universal groupoid $\mathcal{U}_{1}$ of order one. Let $%
\mathcal{U}_{1}{}^{p,q}$ denote the set of all $1$-arrows from $p$ to $q$ so
that $\mathcal{U}_{1}=\cup _{p,q\in M}\mathcal{U}_{1}{}^{p,q}.$ Consider the
group bundle $\mathcal{A}\overset{def}{=}\cup _{p\in M}\mathcal{U}%
_{1}{}^{p,p}$. A choice of coordinates around $p$ identifies $\mathcal{U}%
_{1}{}^{p,p}$ with $GL(n,\mathbb{R)}$. The sections of $\mathcal{%
A\rightarrow }M$ are called gauge transformations and they form a group with
fiberwise composition. A local section $f$ of $\mathcal{A\rightarrow }M$
over $(U,x^{i})$ is of the form $(f_{j}^{i}(x))$ and $f$ is smooth if $%
(f_{j}^{i}(x))$ are smooth functions. This definition does not depend on
coordinates. To be consistent with our notation above, we denote the group
of gauge transformations $\Gamma \mathcal{A}$ by the same letter $\mathcal{A}
$ and always assume that gauge transformations are smooth. It is standard in
gauge theory to construct $\mathcal{A}$ as sections of an appropriate
associated bundle of the principal bundle $\cup _{x\in M}\mathcal{U}%
_{1}{}^{p,x}.$ This construction, apparently more complicated than ours,
however applies to all principal bundles. Note that a coordinate change $%
(x^{i})\rightarrow (y^{i})$ transforms $(f_{j}^{i}(x))$ by $\frac{\partial
y^{i}}{\partial x^{a}}f_{b}^{a}(x)\frac{\partial x^{b}}{\partial y^{j}}%
=f_{j}^{i}(y).$ Therefore $f\in \mathcal{A}$ is an invertible section of $%
T^{\ast }\otimes T\rightarrow M.$

Now let $(M,\varepsilon )$ be a parallelizable manifold with the splitting $%
\varepsilon .$ Some $f\in \mathcal{A}$ acts on $\varepsilon $ as $%
(f\varepsilon )(p,q)\overset{def}{=}f(p)\circ \varepsilon (p,q)\circ
f(q)^{-1}$ where $\circ $ denotes composition of $1$-arrows. In coordinates,
this action is given by

\begin{equation}
(f\varepsilon )_{j}^{i}(x,y)=f_{a}^{i}(y)\varepsilon
_{b}^{a}(x,y)g_{j}^{b}(x)\text{ \ \ }f_{a}^{i}(z)g_{j}^{a}(z)=\delta _{j}^{i}
\end{equation}%
Clearly $f\varepsilon $ is another splitting and therefore $(M,f\varepsilon
) $ is another parallelizable manifold. Let $\Xi $ denote the set of all
splittings on $M.$

\begin{lemma}
$\mathcal{A}$ acts transitively on $\Xi $
\end{lemma}

To prove the Lemma 11, let $\varepsilon _{0},$ $\varepsilon _{1}$ two
splittings. We fix some $p\in M.$ Clearly, for $x\in M$ there exists a
unique $f(x)\in \mathcal{U}^{x,x}$ satisfying $\varepsilon
_{0}(p,x)=f(x)\circ \varepsilon _{1}(p,x)$ which defines some $f\in \mathcal{%
A}.$ Now $\varepsilon _{0}(x,y)=\varepsilon _{0}(p,y)\circ \varepsilon
_{0}(x,p)=g(y)\circ \varepsilon _{1}(p,y)\circ \varepsilon _{1}(x,p)\circ
g(x)^{-1}$ $=g(y)\circ \varepsilon _{1}(x,y)\circ g(x)^{-1}$, that is, $%
\varepsilon _{0}=f\varepsilon _{1}.$

Note, however, that $g$ is by no means unique in $\varepsilon
_{0}=g\varepsilon _{1}.$

Now let $f_{t},$ $0\leq t\leq \epsilon $, be a smooth curve of gauge
transformations starting from the identity, which means that we require the
components $f_{j}^{i}(t,x),$ $x\in (U,x^{i})$ of $f_{t}$ to be smooth
functions and $f_{0}=id_{M}$, that is, $(f_{0})_{j}^{i}(x)=\delta _{j}^{i}.$
If $(M,\varepsilon )$ is a parallelizable manifold and $f_{t}$ is a curve of
gauge transformations, we call $(M,f_{t}\varepsilon )$ a deformation of $%
(M,\varepsilon ).$ It is crucial to observe that a deformation does not
touch the base manifold $M$ but deforms only $1$-arrows (see [6], pg. 183).
Writing $\widetilde{\mathfrak{R}}(\varepsilon ),$ $\widetilde{\nabla }%
(\varepsilon )...$ for $\widetilde{\mathfrak{R}},$ $\widetilde{\nabla }...$
to indicate the dependence on the parallelism$,$ clearly $\widetilde{%
\mathfrak{R}}(f_{t}\varepsilon )=0$ for all $t$ since $f_{t}\varepsilon $
defines a splitting for all $t$ (this can be checked directly using (42)
below).

\begin{definition}
Let $(M,\varepsilon )$ be a local Lie group and $f_{t}$ be a curve. If $%
(M,f_{t}\varepsilon )$ is a local Lie group for all $t,$ then we call $%
(M,f_{t}\varepsilon )$ a deformation of the local Lie group $(M,\varepsilon
).$
\end{definition}

Therefore $(M,f_{t}\varepsilon )$ is a deformation of the local Lie group $%
(M,\varepsilon )\Leftrightarrow \mathcal{R}(f_{t}\varepsilon )=0$ for all $%
t\Leftrightarrow \widehat{\mathfrak{R}}(f_{t}\varepsilon )=0$ for all $%
t\Leftrightarrow \widetilde{\nabla }(f_{t}\varepsilon )T(f_{t}\varepsilon
)=0 $ for all $t.$ To simplify our notation, henceforth we will simply write
$\mathcal{R}_{t},$ $\widetilde{\nabla }_{t}...$ for $\mathcal{R}%
(f_{t}\varepsilon )$, $\widetilde{\nabla }(f_{t}\varepsilon )...$but avoid
the notation $M_{t}$ of [6] for $(M,f_{t}\varepsilon )$ as it gives the
impression that the base manifold $M$ is deformed. We will call $%
(M,f_{t}\varepsilon )$ simply a deformation which we henceforth assume$.$

So if $(M,f_{t}\varepsilon )$ is a deformation, we have the scenario of the
Sections 2, 3, 4 for all $t\geq 0.$ In particular we have the resolution
\begin{equation}
0\rightarrow \widehat{\Theta }_{t}\rightarrow T\overset{\widehat{d}_{t}}{%
\longrightarrow }T^{\ast }\otimes T\overset{\widehat{d}_{t}}{\longrightarrow
}\Lambda ^{2}(T^{\ast })\otimes T\overset{\widehat{d}_{t}}{\longrightarrow }%
....\overset{\widehat{d}_{t}}{\longrightarrow }\Lambda ^{n}(T^{\ast
})\otimes T
\end{equation}%
which reduces to (24) for $t=0$ and restricts to $f_{t}\varepsilon $%
-invariant forms as in (30). We will denote $f_{t}\varepsilon $-invariance
by $\widetilde{\cdot }^{t}$ so that $\widetilde{\cdot }^{0}$ is the same as $%
\widetilde{\cdot }$ and $\widetilde{T_{s}^{r}}^{t}$ denotes $%
f_{t}\varepsilon $-invariant tensor fields. By Proposition 1 $\xi \in
\widetilde{T_{s}^{r}}^{t}\Leftrightarrow \widetilde{\nabla }_{t}\xi =0.$

Note that a deformation deforms only the operators but not the spaces in
(24). Therefore the solution sheaves $\widehat{\Theta }_{t}$ and the
cohomology groups $H^{\ast }(\widetilde{\Theta }_{t},\widetilde{\Theta }%
_{t}),$ which localize at any $p\in M$ to $H^{\ast }((\widetilde{\Theta }%
_{t})_{p},(\widetilde{\Theta }_{t})_{p},$ are also deformed. Clearly, a
deformation imposes strong restrictions on the curve $f_{t}$ and one such
condition will be given in Proposition 14 below.

Now $f_{t}(p)$ defines an isomorphism $f_{t}(p)_{\ast }:T_{p}\rightarrow
T_{p}$ defined by $(\xi _{p}^{i})\rightarrow f_{t}(p)_{a}^{i}\xi _{p}^{a}$
which extends to an isomorphism on tensor space $(T_{p})_{s}^{r}$ and
therefore to tensor fields $T_{s}^{r}.$ In particular, we obtain the diagram

\begin{equation}
\begin{array}{ccccccccc}
\widehat{\Theta }_{t} & \overset{\widehat{d}_{t}}{\rightarrow } & T &
\overset{\widehat{d}_{t}}{\rightarrow } & T^{\ast }\otimes T & \overset{%
\widehat{d}_{t}}{\rightarrow } & ... & \overset{\widehat{d}_{t}}{\rightarrow
} & \Lambda ^{n}(T^{\ast })\otimes T \\
&  & \uparrow _{(f_{t})_{\ast }} &  & \uparrow _{(f_{t})_{\ast }} &  &
\uparrow _{(f_{t})_{\ast }} &  & \uparrow _{(f_{t})_{\ast }} \\
\widehat{\Theta }_{0} & \overset{\widehat{d}_{0}}{\rightarrow } & T &
\overset{\widehat{d}_{0}}{\rightarrow } & T^{\ast }\otimes T & \overset{%
\widehat{d}_{0}}{\rightarrow } & ... & \overset{\widehat{d}_{0}}{\rightarrow
} & \Lambda ^{n}(T^{\ast })\otimes T%
\end{array}%
\end{equation}

We will show the commutatitivity of (40) for constant deformations below
(Proposition 14) but do not know any conditions on $f_{t}$ which makes (40)
commute in general. Note that commutativity of (40) implies that $%
(f_{t})_{\ast }$ induces isomorphism in cohomology and therefore $\dim H^{k}(%
\widehat{\Theta }_{t},\widehat{\Theta }_{t})$ is independant of $t$ (which
is one of the main technical difficulties in [6], see Theorems 4.1, 4.2, 4.3
in [6]).

Reasoning like Kodaira-Spencer, we now ask the following fundamental
questions.

\textbf{Q1)} How to define the derivative (with respect to time $t)$ of a
deformation?

\textbf{Q2) }How to define a constant deformation so that the derivative of
a deformation vanishes for all $t$ if and only if the deformation is
constant.

Even though the main idea of deformation of a complex structure is very
intuitive and simple as explained already in the introduction of [6], the
execution of this idea introduces subtle technical difficulties. For
instance, in addition to the one mentioned above, it turns out that "if and
only if" in \textbf{Q2} is too much to expect. On the other hand, we will
show below (Proposition 17) that \textbf{Q1, Q2} have surprisingly simple
and complete answers in the present framework of local Lie groups.

Now if $f_{t}$ is a curve, then $f_{t}(p)$ is a curve in the Lie group $%
\mathcal{U}^{p,p}\simeq GL(n,\mathbb{R})$ starting from the identity (not
necessarily a $1$-parameter subgroup!!). Therefore $\left[ \frac{df_{t}(p)}{%
dt}\right] _{t=0}\in \mathcal{L(}\mathcal{U}^{p,p})=$ the Lie algebra of $%
\mathcal{U}^{p,p}$. Consider the Lie algebra bundle $\mathfrak{A}\overset{def%
}{=}\cup _{p\in M}\mathcal{L(}\mathcal{U}^{p,p})\rightarrow M$ (actually
another associated bundle) which can be canonically identified with $T^{\ast
}\otimes T\rightarrow M$. If $f\in \mathcal{A},$ then $\frac{df_{t}}{dt}\in
\mathfrak{A}$ and we can view $\mathfrak{A}$ as the Lie algebra of the group
$\mathcal{A}$ even though $\mathcal{A}$ is not a Lie group. If $%
(M,f_{t}\varepsilon )$ is a deformation, it is natural to expect that its
derivative $\frac{df_{t}}{dt}\in \mathfrak{A}$ will be quite relevant for
\textbf{Q1 }and \textbf{Q2. }

It is also worthwhile to note that in [6] we "see" the deformation as an
"outcome" but there is no object that causes this deformation. The presence
of this concrete object $f_{t}$ in our framework will greatly simplify the
picture as already hinted by (40).

\section{Derivative of a deformation}

Let $(M,f_{t}\varepsilon )$ be a deformation. We rewrite (38) for all $t$ as

\begin{equation}
(f_{t}\varepsilon )_{j}^{i}(x,y)=f_{a}^{i}(t,y)\varepsilon
_{b}^{a}(x,y)g_{j}^{b}(t,x)\text{ \ \ }f_{a}^{i}(t,z)g_{j}^{a}(t,z)=\delta
_{j}^{i}
\end{equation}

Now (2) and (41) give

\begin{eqnarray}
\Gamma _{kj}^{i}(f_{t}\varepsilon ) &=&\left[ \frac{\partial }{\partial y^{k}%
}f_{a}^{i}(t,y)\varepsilon _{b}^{a}(x,y)g_{j}^{b}(t,x)\right] _{y=x}  \notag
\\
&=&\frac{\partial f_{a}^{i}(t,x)}{\partial x^{k}}%
g_{j}^{a}(t,x)+f_{b}^{i}(t,x)\Gamma _{kc}^{b}(x)g_{j}^{c}(t,x)
\end{eqnarray}%
For simplicity of notation, we will denote $\Gamma
_{kj}^{i}(f_{t}\varepsilon )$ by $(\Gamma _{t})_{kj}^{i}$ and sometimes drop
the variables $t$ and/or $x$ from our notation in our local computations.
Thus we rewrite the important formula (42) as

\begin{equation}
\left( \Gamma _{t}\right) _{kj}^{i}=\frac{\partial f_{a}^{i}}{\partial x^{k}}%
g_{j}^{a}+f_{b}^{i}\left( \Gamma _{0}\right) _{kc}^{b}g_{j}^{c},\text{ \ }%
t\geq 0
\end{equation}

Having (43) at our disposal, we are now in position to compute everything in
this note explicitly in coordinates including (40). Such computations do not
require any ingeniuity but a lot of patience and may not be much rewarding
at the end.

Now recall the isomorphism $f_{t}(p)_{\ast }:T_{p}\rightarrow T_{p}$ of
vector spaces and the isomorphism $(T_{p},[\cdot ,\cdot ]_{t,p})\simeq $ $(%
\widetilde{\Theta }_{t})_{p}$ of Lie algebras where $[\cdot ,\cdot ]_{t,p}$
is defined by (6) with $T$ replaced by $T_{t}.$ Clearly the vector space
isomorphism $f_{t}(p)_{\ast }:(T_{p},[\cdot ,\cdot ]_{0,p})\rightarrow
(T_{p},[\cdot ,\cdot ]_{t,p})$ need not be a Lie algebra isomorphism.
Equivalently, the Lie algebras $\left\{ \widetilde{\Theta }_{t},\text{ }%
t\geq 0\right\} $ need not be isomorphic (by $(f_{t})_{\ast }!)$ to the
original Lie algebra $\widetilde{\Theta }_{0}$ during the deformation. If
they are, obviously our deformation is not really deforming $\widetilde{%
\Theta }_{0}$.

\begin{definition}
The deformation $(M,f_{t}\varepsilon )$ is constant if $(f_{t})_{\ast }:%
\widetilde{\Theta }_{0}\rightarrow \widetilde{\Theta }_{t}$ is an
isomorphism of Lie algebras for all $t.$
\end{definition}

Note that \textit{any }linear map $l:T_{p}\rightarrow T_{p}$ extends to a
linear map $l_{\ast }:\widetilde{\Theta }_{0}\rightarrow \widetilde{\Theta }%
_{t}$ \ for all $t,$ because the vector fields in the solution sheaves $%
\widetilde{\Theta }_{t}$ are determined by their values at any point by (6).
Therefore, $(M,f_{t}\varepsilon )$ is constant $\Leftrightarrow \lbrack \xi
_{p},\eta _{p}]_{t,p}=[\xi _{p},\eta _{p}]_{0,p}$ for all $\xi _{p},\eta
_{p}\in T_{p}$ for some (hence all) $p\in M.$ In view of (6), this condition
is equivalent to

\begin{equation}
f_{a}^{i}(t,x)\left( T_{0}\right) _{jk}^{a}(x)=\left( T_{t}\right)
_{ab}^{i}(t,x)f_{j}^{a}(t,x)f_{k}^{b}(t,x)\text{ \ }t\geq 0
\end{equation}

(44) simply states that $(f_{t})_{\ast }$ maps the original torsion $T=T_{0}$
to the torsion $T_{t}$ at time $t.$ It is now crucial to observe that (44)
is \textit{not }obtained from (43) by alternating $k,j$ for otherwise all
deformations would be constant. Therefore, it is quite natural to define the
curvature $\kappa (f_{t})$ of $f_{t}$ as a measure of how much the equality
in (44) fails to hold. So we fix some $p\in M$, choose $\xi _{p},\eta
_{p}\in T_{p}$ and define $\kappa (f_{t})(p)(\xi _{p},\eta _{p})$ by

\begin{equation}
\kappa (f_{t})(p)(\xi _{p},\eta _{p})\overset{def}{=}[\xi _{p},\eta
_{p}]_{t,p}-[\xi _{p},\eta _{p}]_{0,p}
\end{equation}

Clearly $\kappa (f_{t})(p)\in \Lambda ^{2}(T_{p}^{\ast })\otimes T_{p},$
that is, $\kappa (f_{t})\in \Lambda ^{2}(T^{\ast })\otimes T$ and $\kappa
(f_{0})=\kappa (id)=0.$ If $\kappa (f_{t})(p)=0$ for some $p,$ then clearly $%
\kappa (f_{t})(p)=0$ for all $p.$

This is a good place to state

\begin{proposition}
If $(M,f_{t}\varepsilon )$ is constant, the following diagram commutes
\end{proposition}

\begin{equation}
\begin{array}{ccccccccc}
& \widetilde{\Theta }_{t}= & \widetilde{T}^{t} & \overset{\widehat{d}_{t}}{%
\rightarrow } & \widetilde{T^{\ast }\otimes T}^{t} & \overset{\widehat{d}_{t}%
}{\rightarrow } & ... & \overset{\widehat{d}_{t}}{\rightarrow } & \widetilde{%
\Lambda ^{n}(T^{\ast })\otimes T}^{t} \\
&  & \uparrow _{(f_{t})_{\ast }} &  & \uparrow _{(f_{t})_{\ast }} &  &
\uparrow _{(f_{t})_{\ast }} &  & \uparrow _{(f_{t})_{\ast }} \\
& \widetilde{\Theta }_{0}= & \widetilde{T}^{0} & \overset{\widehat{d}_{0}}{%
\rightarrow } & \widetilde{T^{\ast }\otimes T}^{0} & \overset{\widehat{d}_{0}%
}{\rightarrow } & ... & \overset{\widehat{d}_{0}}{\rightarrow } & \widetilde{%
\Lambda ^{n}(T^{\ast })\otimes T}^{0}%
\end{array}%
\end{equation}

To prove Proposition 14, we write the action (41) as $(\varepsilon
_{t})(p,q)\circ f_{t}(p)=f_{t}(q)\circ \varepsilon (p,q)$ for some fixed $t.$
Therefore if $\xi \in (T_{p})_{s}^{r}$, then the composions $\varepsilon
_{t}(p,q)_{\ast }$ $\circ (f_{t}(p))_{\ast }$ and $f_{t}(q)_{\ast }\circ
\varepsilon (p,q)_{\ast }$ give the same maps $(T_{p})_{s}^{r}\rightarrow
(T_{q})_{s}^{r}.$ Now it easily follows from definitions that $f_{t}$
induces an isomorphism $(f_{t})_{\ast }:\widetilde{T_{s}^{r}}^{0}\rightarrow
\widetilde{T_{s}^{r}}^{t}$ so that the vertical arrows in (46) are correct.
By Proposition 1, $\xi \in \widetilde{\Lambda ^{k}(T^{\ast })\otimes T}%
^{t}\Rightarrow \lbrack \widetilde{\nabla }_{t}]\xi =0$ so that (28) holds
for all $t.$ Since $(f_{t})_{\ast }:$ $\widetilde{\Lambda ^{k}(T^{\ast
})\otimes T}^{0}\rightarrow \widetilde{\Lambda ^{k}(T^{\ast })\otimes T}^{t}$
and also maps $T_{0}$ to $T_{t}$ by our assumption, the conclusion follows.

Now if we fix $p$ and regard $\kappa (f_{t})(p)$ as a function of $t$, then $%
\kappa (f_{t})(p)=0$ for all $t\Leftrightarrow $ the deformation $%
(M,f_{t}\varepsilon )$ is constant. Clearly $\frac{d}{dt}\kappa
(f_{t})(p)\in \Lambda ^{2}(T_{p}^{\ast })\otimes T_{p}$ and therefore $\frac{%
d}{dt}\kappa (f_{t})\in \Lambda ^{2}(T^{\ast })\otimes T.$

\begin{lemma}
$\frac{d}{dt}\kappa (f_{t})\in \widetilde{\Lambda ^{2}(T^{\ast })\otimes T}%
^{t}$
\end{lemma}

To prove the Lemma, we need to check $\widetilde{\nabla }_{t}(\frac{d}{dt}%
\kappa (f_{t}))=0.$ This is a somewhat long but straightforward computation
using the various definitions involved.

Lemma 15 shows that $\frac{d}{dt}\kappa (f_{t})$ is determined by its value $%
\frac{d}{dt}\kappa _{t}(p)$ at any point $p.$ Now we use the notation $\frac{%
d}{dt}(M,f_{t}\varepsilon )$ for the derivative of the deformation $%
(M,f_{t}\varepsilon )$ at time $t,$ a concept yet to be defined.

\begin{definition}
$\frac{d}{dt}(M,f_{t}\varepsilon )\overset{def}{=}\frac{d}{dt}\kappa (f_{t})$
\end{definition}

Note that our notation in Definition 16 suggests the "equivalence" of the
notations $(M,f_{t}\varepsilon )=\kappa (f_{t})$ which is quite reasonable:
deforming a Lie algebra is the same as deforming its bracket.

As the next proposition shows, Definition 16 trivializes \textbf{Q1, Q2.}

\begin{proposition}
The following are equivalent
\end{proposition}

$i)$ The deformation $(M,f_{t}\varepsilon )$ is constant

$ii)$ $\kappa (f_{t})=0$ for all $t$

$iii)$ $\kappa (f_{t})(p)=0$ for all $t$ and some $p\in M$

$iv)$ $\frac{d}{dt}\kappa (f_{t})=0$ for all $t$

$v)$ $\frac{d}{dt}\kappa (f_{t})(p)=0$ for all $t$ and some $p\in M$

Since $\kappa (f_{t})(p)$ as a function of $t$ takes values in a fixed
vector space, observe the amusing fact that the implication $v)\Rightarrow
iii)$ is implied by calculus in the same way as \textbf{Q2} is inspired by
calculus.

Surprisingly, the trivial Proposition 17 does not even make use of (43) and
at first sight it is unrelated to the complex (30) in contrast to [6] where $%
\frac{d}{dt}M_{t}$ is defined as a cohomology class in $H^{1}$ from the
outset (see pages 188-192). In this direction, we will prove

\begin{proposition}
Let $(M,f_{t}\varepsilon )$ be a constant deformation. Then $\widehat{d}%
_{0}\left( \frac{df_{t}}{dt}\right) _{t=0}=0.$ Therefore $\left( \frac{df_{t}%
}{dt}\right) _{t=0}$ defines a cohomology class in $H^{1}(\Theta _{0},\Theta
_{0}).$
\end{proposition}

Before we start the proof, observe that Proposition 18 is quite
dissappointing for someone whose is familiar with the role of $H^{1}(\Phi )$
in Kodaira-Spencer theory where $\Phi $ is the sheaf of holomorphic vector
fields on a complex manifold. The reason is that $H^{1}(\Phi )$ incorporates
honest deformations and there is a very good reason (see Remark 3 below) to
believe that the map constructed in Proposition 18 (which we will call the
Kodaira-Spencer map) will be surjective, that is, all cohomology classes in $%
H^{1}(\Theta _{0},\Theta _{0})$ arise from constant deformations. Therefore $%
H^{1}(\Theta _{0},\Theta _{0})$ is not much an interesting object. We will
comment on this point further in Remark 3 below.

To prepare for the proof, we ask the conditions imposed on $f_{t}$ by a
deformation $(M,f_{t}\varepsilon )$. Recall that $\mathcal{R}_{t}=0$ holds
identically for all $t$ since we assume that $(M,f_{t}\varepsilon )$ is a
local Lie group for all $t.$ To compute $\mathcal{R}_{t}$, we check the
integrability conditions of

\begin{equation}
\frac{\partial y^{i}}{\partial x^{j}}=f_{a}^{i}(y)\varepsilon
_{b}^{a}(x,y)f_{j}^{b}(x)
\end{equation}%
by differentiating (47) with respect to $x^{k},$ substituting back from (47)
and alternating $k,j.$ After this computation, we differentiate the
expression of $\mathcal{R}_{t}$ with respect to $t$ and set $y=x,$ that is,
we compute $\left[ \frac{d}{dt}(\mathcal{R}_{t})_{kj}^{i}(x,y)\right]
_{y=x}. $ The outcome of this rather long computation is very much
rewarding: all terms cancel except one and what we get is

\begin{equation}
0=\left[ \frac{d}{dt}(\mathcal{R}_{t})_{kj}^{i}(x,y)\right] _{y=x}=\frac{%
df_{a}^{i}(t,x)}{dt}T_{kj}^{a}(x)\text{ \ \ }t\geq 0
\end{equation}

The crucial fact in (48) is that $T$ is the torsion of the original Lie
algebra $\widetilde{\Theta }$ and does not depend on $t.$ To write (48) in a
coordinate free form, let $[\widetilde{\Theta },\widetilde{\Theta }]\simeq $
$[\widetilde{\Theta }_{p},\widetilde{\Theta }_{p}]$ denote the derived
algebra of $\widetilde{\Theta }\simeq \widetilde{\Theta }_{p}$. Now (6) and
(48) show that the linear map $\frac{df_{t}(p)}{dt}\in T_{p}^{\ast }\otimes
T_{p}$ vanishes on $[\widetilde{\Theta }_{p},\widetilde{\Theta }_{p}]$ for
all $t.$ We single out this important fact as

\begin{proposition}
Let $(M,f_{t}\varepsilon )$ be any deformation. Then the linear map $\frac{%
df_{t}(p)}{dt}\in T_{p}^{\ast }\otimes T_{p}$ vanishes on the derived
algebra $[\widetilde{\Theta }_{p},\widetilde{\Theta }_{p}]$ of the original
Lie algebra for all $t\geq 0.$
\end{proposition}

\begin{corollary}
If $[\widetilde{\Theta }_{p},\widetilde{\Theta }_{p}]=\widetilde{\Theta }%
_{p} $ (for instance if $\widetilde{\Theta }_{p}$ is semisimple), then $%
f_{t}(p)=id$ for all $t.$
\end{corollary}

Indeed $[\widetilde{\Theta }_{p},\widetilde{\Theta }_{p}]=\widetilde{\Theta }%
_{p}$ implies $\frac{df_{t}(p)}{dt}=0$ for all $t.$ Proposition 19 shows
that the smaller the derived algebra, the more freedom we have in deforming
the Lie algebra. In particular, semisimple Lie algebras resist deformations
(even constant ones other than the identity! This peculiarity occurs because
we are deforming by gauge transformations, see Remark 1).

Now suppose that the deformation is constant. In (44) we substitute $T_{t}$
from (43), simplify, and omitting $t,x$, we rewrite (44) as

\begin{equation}
f_{a}^{i}T_{jk}^{a}=\left[ \frac{\partial f_{k}^{i}}{\partial x^{a}}%
f_{j}^{a}+f_{c}^{i}\Gamma _{ak}^{c}f_{j}^{a}\right] _{[kj]}
\end{equation}

Now differentiation of (49) with respect to $t$ at $t=0$ gives

\begin{equation}
\frac{df_{a}^{i}}{dt}T_{jk}^{a}=\left[ \frac{\partial }{\partial x^{a}}%
\left( \frac{df_{k}^{i}}{dt}\right) f_{j}^{a}+\frac{\partial f_{k}^{i}}{%
\partial x^{a}}\frac{df_{j}^{a}}{dt}+\frac{df_{c}^{i}}{dt}\Gamma
_{ak}^{c}f_{j}^{a}+f_{c}^{i}\Gamma _{ak}^{c}\frac{df_{j}^{a}}{dt}\right]
_{[kj],\text{ }t=0}
\end{equation}

The LHS of (50) vanishes by (48). Similarly, since $f_{j}^{a}=\delta
_{j}^{a} $ for $t=0,$ the third term inside the paranthesis vanishes upon
alternation. Since $f_{k}^{i}=\delta _{k}^{i}$ for all $x$ at $t=0,$ we $%
\frac{\partial f_{k}^{i}}{\partial x^{a}}=0$ and the second term vanishes
too and (50) becomes

\begin{eqnarray}
0 &=&\left[ \frac{\partial }{\partial x^{j}}\left( \frac{df_{k}^{i}}{dt}%
\right) +\Gamma _{ak}^{i}\frac{df_{j}^{a}}{dt}\right] _{[jk],\text{ }t=0}
\notag \\
&=&\left[ \widehat{d}_{j}\left( \frac{df_{k}^{i}}{dt}\right) \right] _{[jk],%
\text{ }t=0}  \notag \\
&=&\widehat{d}_{0}\left( \frac{df}{dt}\right) _{t=0}
\end{eqnarray}%
proving Proposition 18.

In the above proof, we actually use only $\left( \frac{d}{dt}\kappa
_{t}(p)\right) _{t=0}=0.$ Therefore, it is natural to consider the
coefficients $\left( \frac{d^{k}}{dt^{k}}\kappa _{t}(p)\right) _{t=0}$ $\in
\Lambda ^{2}(T^{\ast })\otimes T$ in the Taylor expansion of $\kappa _{t}(p)$
at $t=0$ as an approximation to $\kappa _{t}(p)$. Clearly all these
coefficients vanish for constant deformations which suggests their relevance
for the Remark 3 below.

We will conclude this note with three remarks.

1) In [5], we defined the symmetry group $\mathcal{S}(\mathcal{G)}$ of a
local Lie group $(M,\varepsilon )$ which is characterized by the following
property: $\mathcal{S}(\mathcal{G)}$ is a second order pseudogroup on $M$
which contains the pseudogroup $\mathcal{G}$ as a normal subpseudogroup (if $%
M$ is compact and simply connected, then both $\mathcal{S}(\mathcal{G)}$ and
$\mathcal{G}$ are Lie groups and $\mathcal{G\vartriangleleft S}(\mathcal{G))}
$ and is maximal with respect to this property, that is, any pseudogroup $%
\mathcal{S}$ on $M$ with $\mathcal{G\vartriangleleft S}$ is contained in $%
\mathcal{S}(\mathcal{G)}$. The Lie algebra of vector fields of $\mathcal{S}(%
\mathcal{G)}$ are solutions of a second order linear $PDE$ and is isomorphic
to the semi-direct product $Der(\widetilde{\Theta }_{p})\times \widetilde{%
\Theta }_{p}.$ Let $\gamma (t)$ be a $1$-parameter subgroup in $\mathcal{S}(%
\mathcal{G)}$ (assume for simplicity that $M$ is compact and simply
connected). Now $\gamma (t)$ deforms $\mathcal{G}$ (and therefore $%
\widetilde{\Theta })$ by conjugation as $\mathcal{G\rightarrow \gamma (}t)%
\mathcal{G\gamma }(t)^{-1}.$ The advantage of these deformations is that
they immediately give a surjective map to $H^{1}(\widetilde{\Theta }_{p},%
\widetilde{\Theta }_{p})$ with the required kernel thus completely
clarifying the geometric meaning of $H^{1}(\widetilde{\Theta }_{p},%
\widetilde{\Theta }_{p}).$ The disadvantage is that we can produce only
constant deformations in this way. Therefore, we believe that the
Kodaira-Spencer map in Proposition 18 is also surjective. We believe that
problems related to $H^{1}(\widetilde{\Theta }_{p},\widetilde{\Theta }_{p})$
are important only for understanding the mechanism behind deformations and
developing an invariant formalism to attack more complicated problems but
not for their own sake.

2) As we remarked above, for a curve $f_{t}$ of gauge transformations, the
curves $f_{t}(p)$ need not be $1$-parameter subgroups. Let us call $f_{t}$
regular if this condition is satisfied at all points. We will make the
following conjecture

\textbf{C:} Let $(M,f_{t}\varepsilon )$ be a deformation and $M$ compact.
Then for any sufficiently small $t_{0},$ there exists a deformation $%
(M,g_{s}\varepsilon )$ such that $g_{s}$ is regular and $g_{s_{0}}%
\varepsilon =f_{t_{0}}\varepsilon $ for some $s_{0}.$

Shortly, \textbf{C }states that it is sufficient to work with regular
deformations. Observe that we have $f_{t}\left( \frac{df_{t}}{dt}\right)
_{t=0}=\frac{df_{t}}{dt}$ for regular deformations. Also observe the
similarity of the statement of \textbf{C} to the local surjectivity of the
exponential map on a Lie group. We believe that the truth of \textbf{C} will
allow us to use the full power of the exponential map and greatly simplify
the deformation theory of local Lie groups.

3) The integration theorem (first order systems with initial conditions)
which we repeatedly use here boils down to the well known Frobenious theorem
which holds also in the complex category. Therefore everything in this note
applies to complex local Lie groups and therefore complex Lie groups. So let
$(M,\varepsilon )$ be a complex local Lie group and $\Phi $ be the sheaf of
holomorphic vector fields on $M$ which contains $\widehat{\Theta }$ as a
subsheaf. The exact sequence of sheaves $0\rightarrow \widehat{\Theta }%
\rightarrow \Phi \rightarrow \Phi /\widehat{\Theta }\rightarrow 0$ gives the
long exact sequence of cohomology groups

\begin{equation}
...\rightarrow H^{0}(\Phi /\widehat{\Theta })\rightarrow H^{1}(\widehat{%
\Theta })=H^{1}(\widetilde{\Theta },\widetilde{\Theta })\rightarrow
H^{1}(\Phi )\rightarrow H^{1}(\Phi /\widehat{\Theta }))\rightarrow H^{2}(%
\widehat{\Theta }))\rightarrow
\end{equation}

Now Kodaira-Spencer construct cohomology classes in $H^{1}(\Phi ).$ The
dimension $e$ of the moduli space of complex structures (whenever defined)
satisfies $e\leq \dim H^{1}(\Phi )$ and the driving force of the
Kodaira-Spencer adventure is to prove $e=\dim H^{1}(\Phi )$ which eventually
turns out to be false (see [6], pg. 319). From our standpoint, the
interesting question is whether the analytic structure of a complex local
Lie group $\mathcal{G}$ determines the isomorphism class of its Lie algebra $%
\widetilde{\Theta }$. Now (52) suggests that the role of $H^{1}(\Phi )$ for
complex manifolds is played by $H^{2}(\widetilde{\Theta })$ for local Lie
groups. So the answer should be negative if we can relate nonconstant
deformations to $H^{2}(\widehat{\Theta }).$ Note that up to now we always
assumed the existence of a deformation and derived some consequences. But
how do we construct deformations and what is our freedom in constructing
them? These questions will be studied in some future work.

\bigskip

\textbf{References}

\bigskip

[1] E.Abado\u{g}lu, E.Orta\c{c}gil: Intrinsic characteristic classes of a
local Lie group, Portugal. Math. (N.S), Vol.67, Fasc. 4, 2010, 453-483

[2] C.Chevalley, S.Eilenberg: Cohomology theory of Lie groups and Lie
algebras, Trans. Amer. Math. Soc. 63 (1948), 85-124

[3] R.L.Fernandes: Lie algebroids, holonomy and characteristic classes, Adv.
Math. 170, (2002), 119-179

[4] M.Gerstenhaber: The cohomology structure of an associative ring, Ann. of
Math. (2) 78 (1963), 267-288

[5] G.Karaali, P.J.Olver, E.Orta\c{c}gil, M.A.Ta\c{s}k\i n: The symmetry
group of a local Lie group, in progress

[6] K.Kodaira: Complex manifolds and deformation of complex structures,
Classics in Mathematics, Springer, 1986

[7] K.Kodaira, D.C.Spencer: On deformations of complex analytic structures,
I-II, III, Annals of Math. 67 (1958), 328-466; 71. (1960), 43-76

[8] K.McKenzie: Lie groupoids and Lie algebroids in differential geometry,
London Math. Soc. Lecture Note Series, 124, Cambridge University Press, 1987

[9] P.J.Olver: Non-associative local Lie groups, J. Lie Theory 6 (1996),
23-51

[10] E.Orta\c{c}gil: Riemannian geometry as a curved prehomogeneous
geometry, arXiv: 1003.322

[11] J.F.Pommaret: Systems of partial differential equations and Lie
pseudogroups, Gordon \& Breach Science Publishers, 1978

[12] D.C.Spencer: Deformation of structures on manifolds defined by
transitive continuous pseudogroups, Part I: Infinitesimal deformations of
structure, Part II: Deformations of structure, Annals of Math. Vol. 76. No.
2, 1962, 362, 306-445

[13] H.C.Wang: Closed manifolds with homogeneous complex structures, Amer.
J. Math., 76 (1954), 1-32

\bigskip

Erc\"{u}ment Orta\c{c}gil, Bo\u{g}azi\c{c}i University, Math. Dept. (Emer.),
Bebek, \.{I}stanbul, 34342, T\"{u}rkiye

e-mail: ortacgil@boun.edu.tr

\end{document}